\documentclass{amsart}
\usepackage[dvips]{graphicx}
\usepackage{amscd}
\usepackage{amsmath}
\usepackage{amsxtra}
\usepackage{amsfonts}
\usepackage{amssymb}
\newtheorem{theorem}{Theorem}[section]

\newtheorem{lemma}[theorem]{Lemma}

\theoremstyle{definition}
\newtheorem{definition}[theorem]{Definition}
\newtheorem{remark}[theorem]{Remark}

\theoremstyle{remark}

\renewcommand{\theclaim}{\textup{\theclaim}}

\newtheorem*{acknowledgements}{Acknowledgements}

\numberwithin{equation}{section}

\input cyracc.def

\begin{document}

\title[The Laplacian subalgebra of $\mathcal{L}(\mathbb{F}_N)^{\overline{\otimes}_k}$ is a strongly singular MASA]{The Laplacian subalgebra
of $\mathcal{L}(\mathbb{F}_N)^{\overline{\otimes}_k}$ is a strongly singular MASA}
\author{Teodor \c Stefan B\^ ildea}
\address{Department of Mathematics\\
The University of Iowa\\
14 MacLean Hall\\
Iowa City, IA 52242-1419\\
U.S.A.\\} \email{stefan-bildea@uiowa.edu}


\begin{abstract}
Using the techniques of A. Sinclair and R.Smith, we show that the Laplacian (radial) von Neumann subalgebra of
$\mathcal{L}(\mathbb{F}_N)^{\overline{\otimes}_k}$ is a strongly singular maximal abelian subalgebra for every $k\ge 1$.
\end{abstract}

\maketitle

\section{\label{Intro}Introduction} Consider a type ${\rm II}_1$ factor $M$ with trace $\tau$, and let $A\subset M$ be a maximal
abelian von Neumann subalgebra (masa) of $M$. Following Dixmier (\cite{Dix}), let $N_M(A)=\{ u\in A| u\mbox{ unitary },uAu^{\ast}=A\}$ be the
group of normalizing unitaries of $A$.  According to the size of  $B=N_M(A)^{"}$ in $M$, two extreme situations can occur: $B=M$, that is $B$
has enough unitaries to generate $M$ and in this case $A$ is called {\bf regular} or ${\bf Cartan}$;  $B=A$, in which case the only normalizing
unitaries are the unitaries in $A$, and $B$ is called {\bf singular}. Regular/Cartan masas appear in a natural way in cross-products:
$L^{\infty}[0,1]$ is such a Cartan masa in its crossed product with some countable group $G$. On the other hand, until recently, singular masas
where hard to find. Dixmier, who started the investigation of masas, gave the example of the singular masa inside the free group factor
$\mathcal{L}(\mathbb{F}_2)$. The masa is generated by one of the two generators of $\mathbb{F}_2$ inside the factor. Building on work of Pytlik
(\cite{Pytlik}) and Popa (\cite{Popa},\cite{Popa2}), R\u adulescu \cite{FR} showed, using the Pukansky invariant,  that the Laplacian (a.k.a
radial) subalgebra of $\mathcal{L}(\mathbb{F}_N)$ is a singular masa. Together with Boca (\cite{BR}), they gave more examples of singular masas
in group von Neumann algebras of certain free products of groups. It was Popa who actually showed in \cite{Popa} that every separable $II_1$
factor has a singular masa. On the other hand, it is worth mentioning that Voiculescu proved the absence of Cartan masas in free group factors.
\par The class of strongly singular  masas, introduced by Sinclair and Smith in  \cite{ss2} turn out to be more tractable and a plethora of examples
can be found now in \cite{ss},\cite{ss2},\cite{GSS}, mostly in the context of group von Neumann algebras or crossed products with groups.
Furthermore, using the notion of asymptotic homomorphism and good criteria to detect such maps as well as strong singularity in general,
examples where produced demonstrating strong singularity for certain abelian algebras. In particular, among the first examples of strongly
singular masas were the ones known to Dixmier (as singular masas) and the Laplacian subalgebra of the free group factor.
\par The aim of this paper is to present yet another set of examples of strongly singular masas living inside tensor products of
$\mathcal{L}(\mathbb{F}_N)$ with itself. The examples are natural generalizations of the Laplacian subalgebra of $\mathcal{L}(\mathbb{F}_N)$,
generated by the selfadjoint element
\[ w=\sum_{i=1}^N(g_i+g_i^{-1}),\]
where we have identified the generators $g_1,...,g_N$ of $\mathbb{F}_N$ with the unitaries in $\mathcal{L}(\mathbb{F}_N)$ given by the left
regular representation. In $\mathcal{L}(\mathbb{F}_N)^{\overline{\otimes}_k}$ we consider  the abelian von Neumann subalgebra generated by
\[ w^{(k)}=\sum_{i=1}^N(g_i^{\otimes_k}+(g_i^{-1})^{\otimes_k}),\]
and show that the unique trace-preserving conditional expectation onto this subalgebra is an asymptotic isomorphism, and hence the subalgebra
is a strongly singular masa for each $k\ge 1$.\\
The article is organized as follows: after some preliminaries in section \ref{pre}, we prove in section \ref{lemmas} an essential lemma on free groups.
In section \ref{ssm} we review the notion of strongly singular masa and the tools needed to prove our result.  In the last section we present the main
result.

\section{\label{pre}Preliminaries and notation}
We denote by $\mathbb{F}_N$ the free group on a fixed set of $N\ge 2$ generators $g_1,g_2,...,g_N$. We will call
$S_N=\{g_1,g_2,...,g_N,g_1^{-1},g_2^{-2},...,g_N^{-1}\}$ a generating set for $\mathbb{F}_N$. We denote with $e$ the empty word.
On $\mathbb{F}_N$ we will consider the canonical length function $|\cdot|$, defined by $|e|=0$ and $|w|=|\alpha_1|+|\alpha_2|+...+|\alpha_p|$
if $w\in\mathbb{F}_N$ has the reduced form $g_{i_1}^{\alpha_1}g_{i_2}^{\alpha_2}...g_{i_p}^{\alpha_p}(i_j\ne i_{j+1})$.

The von Neumann group algebra $\mathcal{L}(\mathbb{F}_N)$ is the norm closure of the left regular representation of $\mathbb{F}_N$ on
$l^2(\mathbb{F}_N)$. It is well known that $\mathcal{L}(\mathbb{F}_N)$ is a type ${\rm II}_1$ factor acting standardly on $l^2(\mathbb{F}_N)$
and that with this identification the $\|\cdot\|_{\tau}$ norm coincides with the usual norm $\|\cdot\|_2$ on $l^2(\mathbb{F}_N)$. We will use
the same notation $\|\cdot\|_2$ for the norm $\|\cdot\|_{\tau}$ given by $\|x\|_{\tau}=(\tau(x^{\ast}x)^{1/2})$. We will identify elements in
$\mathbb{F}_N$ with their left regular representations. We will also identify the elements of $\mathcal{L}(\mathbb{F}_N)$ with the corresponding
vectors in $L^2(\mathcal{L}(\mathbb{F}_N),\tau)=l^2(\mathbb{F}_N)$.
\par If $H$ is a Hilbert space, we denote by $H^{\otimes_k}$ the $k$-folded tensor product of $H$ with itself (with the standard convention
$H^{\otimes_0}:=\mathbb{C}$, $H^{\otimes_1}:=H$). An element $v\in H$ embeds as $v\otimes ...\otimes v=:v^{\otimes_k}\in H^{\otimes_k}$. For a
von Neumann algebra $M$, $M^{\overline{\otimes}_k}$ denotes the $k$-folded von Neumann algebra tensor product of $M$ with itself
($M^{\overline{\otimes}_0}:=\mathbb{C}$, $M^{\overline{\otimes}_1}:=M$). For $x\in M$,  $x\otimes...\otimes x=:x^{\otimes_k}\in
M^{\overline{\otimes}_k}$.

For $n\ge 0$ we define $w^{(k)}_n$ to be the sum of all $k$-folded tensors of reduced words in $\mathbb{F}_N$ of length $n$:
\[ w^{(k)}_0=e^{\otimes_k},\]
\[ w^{(k)}_1=\sum_{i=1}^Ng^{\otimes_k}_i+\sum_{i=1}^N(g_1^{-1})^{\otimes_k},\]
\[ w^{(k)}_n=\sum_{|v|=n}v^{\otimes_k}.\]
The following relations are well known for $k=1$ and hold for any $k\ge 1$:
\begin{equation}\label{wn}
 w^{(k)}_1w^{(k)}_n=w^{(k)}_nw^{(k)}_1=w^{(k)}_{n+1}+(2N-1)w^{(k)}_{n-1}.
\end{equation}
These relations show that the von Neumann algebra $\mathcal{B}^{(k)}$ generated by $w^{(k)}_1$ contains each $w^{(k)}_n$, and hence is
the weak closure of the span of these elements:
\[ \mathcal{B}^{(k)}=\overline{Sp\{{w^{(k)}_n|n\ge 0}\}}^{ w}\subset \mathcal{L}(\mathbb{F}_N)^{\overline{\otimes}_k}.\]
We call $\mathcal{B}^{(k)}$ the {\bf Laplacian} or {\bf radial subalgebra}, in analogy with the case $k=1$.
 The selfadjoints $w^{(k)}_n$ are pairwise orthogonal with respect to the trace $\tau^{\otimes_k}$, and a simple counting argument shows that :
 \begin{equation}\label{norma wn}
  \|w^{(k)}_n\|^2_2=2N(2N-1)^{n-1},\quad n\ge 1,
  \end{equation}
  which represents the number of words of length $n$.
  By the previous remarks
  \[ \{\frac{w^{(k)}_n}{\|w^{(k)}_n\|_2}|n\ge 0\} \]
   is an orthonormal basis for $L^2(\mathcal{B}^{(k)},\tau)$.
   Moreover, if $\mathbb{E}_k:\mathcal{L}(\mathbb{F}_N)^{\overline{\otimes}_k}\to \mathcal{B}^{(k)}$ is the unique  trace-preserving conditional
   expectation, then
 \begin{equation}\label{exp x}
 \mathbb{E}_k(x)=\sum_{n=0}^{\infty}\tau^{\otimes_k}(xw^{(k)}_n)\frac{w^{(k)}_n}{\|w^{(k)}_n\|^2_2}.
 \end{equation}
 In particular, when $v$ is a reduced word of length $p$:
 \begin{equation}\label{vdekori}
  \mathbb{E}_k(v^{\otimes_k})=\frac{w^{(k)}_p}{\|w^{(k)}_p\|^2_2}.
  \end{equation}

\begin{lemma}\label{nonzero} For $k\ge 1$ let $x_1,...,x_{k}\in \mathbb{F}_N$.
With the above notation ,
\[ \mathbb{E}_k(x_1\otimes...\otimes x_k)\ne 0\Leftrightarrow x_1=...=x_k.\]
\end{lemma}
{\it Proof: }
Recall from (\ref{exp x}) that
\[ \mathbb{E}_k(x_1\otimes...\otimes x_k)= \sum_{n=0}^{\infty}\tau^{\otimes_k}(x_1\otimes...\otimes x_kw^{(k)}_n)\frac{w^{(k)}_n}{\|w^{(k)}_n\|^2_2}.\]
Therefore
$\mathbb{E}_k(x_1\otimes...\otimes x_k)\ne 0$
if and only if there is some $n\ge 0$ so that :
\[\tau^{\otimes_k}(x_1\otimes...\otimes x_kw^{(k)}_n)\ne 0.\]
Because of the definitions of $w_n^{(k)}$ and $\tau^{\otimes_k}$ this can happen if and only if for some $n\ge 0$:
\[  \sum_{|v|=n}\tau^{\otimes_k}(x_1v\otimes...\otimes x_kv)\ne 0,\]
or
\[ \sum_{|v|=n}\tau(x_1v)...\tau(x_kv)\ne 0.\]
But this is true exactly when for some $v$ of length $n\ge 0$ we have
\[ x_1=...=x_k=v^{-1}. \quad \mbox{\rule{5pt}{5pt}}\]

\section{\label{lemmas}Counting words in $\mathbb{F}_N$}
In this section we prove a technical lemma about solutions of a conjugacy equation in $\mathbb{F}_N$. If $a,b\in \mathbb{F}_N$ are two words and
no cancellations occur in their concatenation, we will denote this by $a\cdot b$. In general, the product (concatenation) will be denoted $ab$,
meaning that cancellations may or may not take place.
\begin{lemma}\label{mylemma}
Let $a,b\in \mathbb{F}_N$ be non-trivial words.\\
(i) The equation $x\cdot a=b\cdot x$ has  in $\mathbb{F}_N$ at most one solution of fixed positive length $l\ge 1$.
\\ (ii)The equation $xa=bx$
has  in $\mathbb{F}_N$ at most one solution of fixed positive length $l>|a|+|b|$.
\end{lemma}
{\it Proof:} (i) Fix a length $l\ge 1$.
Suppose that a non-trivial word $x\in \mathbb{F}_N$ satisfies $x\cdot a=b\cdot x$ and has length $|x|=l$. This is possible only if  $|a|=|b|$. \par
If $x$ is shorter than $a,b$ then it is a suffix for $a$ and and a prefix for $b$. Furthermore, if $a=a_1\cdot x$ with $a_1\in \mathbb{F}_N$ then
$b=x\cdot a_1$ and $x=a_1^{-1}a$ is the unique solution. As we can see, the solution of prescribed length exists only for particular $a,b$.
It depends also on the desired length, so if it exists, it is unique by construction.
\par If  the length $l$ of $x$ is greater than $|a|=|b|$, then $x$ starts with $b$. Suppose $x=b^m\cdot x_1$, $m\ge 1$ maximal with this property.
This forces $|x_1|<|b|=|a|$ and $x_1\cdot a=b\cdot x_1$. As shown in the first paragraph, there is at most one solution $x_1$, hence at most one
solution $x$ of prescribed length $l$.\\
(ii) Next we consider possible cancellations, that is assume $x$ of length $l>|a|+|b|$ satisfies $xa=bx$. This means that
 there are words $a_1,a_2,b_1,b_2$ and $x_1$ such that $a=a_1\cdot a_2\ne 1$, $
b=b_2\cdot b_1\ne 1$, and
$x=b^{-1}_1\cdot x_1\cdot a^{-1}_1,\quad |x_1|\ge 1,$
and at least one of $a_1$ or $b_1$ is a nonempty word. Looking at the lengths in the resulting equation
\[ b^{-1}_1\cdot x_1\cdot a_2=b_2\cdot x_1\cdot a^{-1}_1\]
we see that:
\begin{equation}\label{lungimi}
|b^{-1}_1|+|a_2|=|b_2|+|a_1^{-1}|\Leftrightarrow |a_1|-|a_2|=|b_1|-|b_2|.
\end{equation}
These are necessary conditions on $a,b$ for the existence of a solution.
By hypothesis, we have $|a_1|>0$ or $|b_1|>0$. Note also that $|a_2a_1|>0,|b_1b_2|>0$. If $|a_1|<|a_2|$, then (\ref{lungimi}) implies $|b_1|<|b_2|$
and the equation can be written:
\[ x_1\cdot (a_2a_1)=(b_1b_2)\cdot x_1.\]
The case $|a_1|>|a_2|$ leads to
\[ (b_2^{-1}b_1^{-1})\cdot x_1=x_1\cdot (a_1^{-1}a_2^{-1}). \]
As before, the above equations, and hence the original equation, have at most one solution of prescribed length. The case $|a_1|=|a_2|$ is not
possible, since it leads to $a_1^{-1}=a_2$ and we have assumed $a=a_1\cdot a_2$ is non-trivial.
\rule{5pt}{5pt}

\section{\label{ssm}Strongly singular masa's}
A maximal abelian selfadjoint subalgebra $\mathcal{A}$ in a ${\rm II}_1$ factor $\mathcal{M}$ is singular if any unitary $u\in \mathcal{M}$
which normalizes $\mathcal{A}$ (i.e. $u\mathcal{A}u^{\ast}=\mathcal{A}$), must lie in $\mathcal{A}$.

To define the notion of strongly singular masa, consider a linear map $\phi:\mathcal{M}_1\to \mathcal{M}_2$ between two type ${\rm II}_1$ factors.
There are several norms for $\phi$, depending on the norms considered on the two algebras. When $\mathcal{M}_1$ has the operator norm and
$\mathcal{M}_2$ has the $\|\cdot\|_2$-norm given by the trace, we denote the resulting norm for $\phi$ by $\|\phi\|_{\infty,2}$, following \cite{ss2}.
\begin{definition} Suppose $\mathcal{A}$ is a masa  in a type ${\rm II}_1$ factor $\mathcal{M}$.\\
1)$\mathcal{A}$ is called $\alpha$-strongly singular (or simply strongly singular for $\alpha=1$) if
\[ \|\mathbb{E}_{u\mathcal{A}u^{\ast}}-\mathbb{E}_{\mathcal{A}}\|_{\infty,2}\ge \alpha\|u-\mathbb{E}_{\mathcal{A}}(u)\|_2,\]
for all unitaries $u\in \mathcal{M}$.\\
2) The conditional expectation $\mathbb{E}_{\mathcal{A}}$ is an asymptotic homomorphism if there is a unitary $u\in \mathcal{A}$ such that
\[ \lim_{|k|\to \infty}\|\mathbb{E}_{\mathcal{A}}(xu^ky)-\mathbb{E}_{\mathcal{A}}(x)\mathbb{E}_{\mathcal{A}}(y)u^k\|_2=0\]
for all $x,y\in \mathcal{M}$.
\end{definition}

  The following provides the criteria we will use in the last section to prove the main result of this paper.
\begin{theorem}\label{criteria}(Sinclair and Smith) Let $\mathcal{A}$ be a abelian von Neumann subalgebra of a type ${\rm II}_1$ factor
 $(\mathcal{M},tr)$, and suppose that there is a $\ast$-isomorphism $\pi:\mathcal{A}\to L^{\infty}[0,1]$ which induces an isometry from
 $L^2(\mathcal{A},tr)$ onto $L^2[0,1]$. Let $\{v_n\in \mathcal{A}|n\ge 0\}$ be an orthonormal basis for $L^2(\mathcal{A},tr)$, and let
 $Y\subseteq \mathcal{M}$ be a set whose linear span is norm dense in $L^2(\mathcal{M},tr)$. Let $\mathbb{E}_{\mathcal{A}}:\mathcal{M}\to \mathcal{A}$
 be the unique conditional expectation satisfying $tr \circ \mathbb{E}_{\mathcal{A}}=tr$. If
\[ \sum^{\infty}_{n=0}\| \mathbb{E}_{\mathcal{A}}(xv_ny)-\mathbb{E}_{\mathcal{A}}(x)\mathbb{E}_{\mathcal{A}}(y)v_n\|_2^2<\infty\]
for all $x,y\in Y$, then $\mathbb{E}_{\mathcal{A}}$ is an asymptotic homomorphism and $\mathcal{A}$ is a strongly singular masa.
\end{theorem}

 \begin{remark}\label{r1}
  As abelian von Neumann algebras, $\mathcal{B}^{(k)}\subset B(L^2(\mathcal{B}^{(k)},\tau^{\otimes_k}))$ and
  $\mathcal{B}^{(1)}\subset B( L^2(\mathcal{B}^{(1)},\tau))$ are $\ast$-isomorphic and the isomorphism induces an isometry between
  the $L^2$-spaces, sending $w^{(k)}_n$ to $w^{(1)}_n, n\ge 0$. In \cite{ss} it is shown that $\mathcal{B}^{(1)}$ and
  $\mathbb{E}_1:\mathcal{L}(\mathbb{F}_N)\to \mathcal{B}^{(1)}$ satisfy the hypothesis of the above theorem (so $\mathbb{E}_1$
  is an asymptotic homomorphism and $\mathcal{B}^{(1)}$ is a strongly singular masa). Via the $\ast$-isomorphism between $\mathcal{B}^{(k)}$
  and $\mathcal{B}^{(1)}$ it follows that $\mathcal{B}^{(k)}$ also satisfies the hypothesis of the theorem, for all $k\ge 2$.
  \end{remark}
We further exploit the $\ast$-isomorphism between $\mathcal{B}^{(k)}$ and  $\mathcal{B}^{(1)}$ in the following :
\begin{lemma}\label{k0}Let $x,y\in \mathbb{F}_N$ with $|x|=l\ge 2,|y|=m\ge 2$. The following equality holds for every $k\ge 1,n\ge l+m$:
\begin{equation}\label{eq1}
\|\mathbb{E}_k(x^{\otimes_k}w^{(k)}_ny^{\otimes_k})-\mathbb{E}_k(x^{\otimes_k})\mathbb{E}_k(y^{\otimes_k})w^{(k)}_n\|^2_2=\|\mathbb{E}_1(xw^{(1)}_ny)-
\mathbb{E}_1(x)\mathbb{E}_1(y)w^{(1)}_n\|^2_2.
\end{equation}

\end{lemma}
 {\it Proof: } To simplify the notation, let $w^{(1)}_n:=w_n$.
Note that if $M$ is a von Neumann algebra with trace tr and $a,b\in M$, then $\|a-b\|^2_2=\|a\|^2_2+\|b\|^2_2 -2$Re$(tr(a^{\ast}b))$. Hence to prove
(\ref{eq1}) means to show
\begin{eqnarray}
\|\mathbb{E}_k(x^{\otimes_k}w^{(k)}_ny^{\otimes_k})\|^2_2&=&\|\mathbb{E}_1(xw_ny)\|^2_2 \label{eq21}\\
\|\mathbb{E}_k(x^{\otimes_k})\mathbb{E}_k(y^{\otimes_k})w^{(k)}_n\|^2_2&=&\|\mathbb{E}_1(x)\mathbb{E}_1(y)w_n\|^2_2\mbox{ and } \label{eq22}\\
\quad \tau^{\otimes_k}((\mathbb{E}_k(x^{\otimes_k}w^{(k)}_ny^{\otimes_k})^{\ast}\mathbb{E}_k(x)\mathbb{E}_k(y)w^{(k)}_n)&=
&\tau((\mathbb{E}_1(xw_ny)^{\ast}\mathbb{E}_1(x)\mathbb{E}_1(y)w_n). \label{eq23}
\end{eqnarray}
Let $\mu(r,s,n;x,y)$ be the number of reduced words in the product $xw_ny$ which result from $r$ cancellations on the left and $s$
cancellations on the right. Remark that we have assumed $n\ge l+m$.
We have:
\[ \mathbb{E}_k(x^{\otimes_k}w^{(k)}_ny^{\otimes_k})=\sum_{r=0}^l\sum_{s=0}^m\mu(r,s,n;x,y)w^{(k)}_{m+l+n-2(r+s)}\|w^{(k)}_{m+l+n-2(r+s)}\|^{-2}_2.\]
Because the selfadjoint elements $w^{(k)}_j,j\ge 0$ are pairwise orthogonal with respect to the trace, we get:
\[\|\mathbb{E}_k(x^{\otimes_k}w^{(k)}_ny^{\otimes_k})\|^2_2=\sum_{r=0}^l\sum_{s=0}^m\left(\frac{\mu(r,s,n;x,y)}{\|w^{(k)}_{m+l+n-2(r+s)}\|_2}\right)^2\]
and (\ref{eq21}) follows using (\ref{norma wn}). Next
\[ \| \mathbb{E}_k(x^{\otimes_k})\mathbb{E}_k(y^{\otimes_k})w^{(k)}_n\|^2_2=\frac{1}{\|w^{(k)}_l\|^4_2}\cdot \frac{1}{\|w^{(k)}_m\|^4_2}\cdot
\|w^{(k)}_lw^{(k)}_mw^{(k)}_n\|^2_2=\]
\[=\frac{1}{\|w_l\|^4_2\cdot \|w_m\|^4_2}\cdot \sum_{\tiny \begin{array}{c}
                                            |v_1|=|t_1|=n\\
                                            |v_2|=|t_2|=l\\
                                            |v_3|=|t_3|=m
                                        \end{array}}
                                        (\tau(v_1v_2v_3t_3t_2t_1))^k=\frac{1}{\|w_l\|^4_2\cdot \|w_m\|^4_2}\cdot \|w_lw_mw_n\|^2_2=\]
\[=\|\mathbb{E}_1(x)\mathbb{E}_1(y)w_n\|^2_2,\]
and this proves (\ref{eq22}). Finally,

\[\tau^{\otimes_k}((\mathbb{E}_k(x^{\otimes_k}w^{(k)}_ny^{\otimes_k}))^{\ast}\mathbb{E}_k(x)\mathbb{E}_k(y)w^{(k)}_n)=\]

\[=\tau^{\otimes_k}((\mathbb{E}_k(x^{\otimes_k}w^{(k)}_ny^{\otimes_k}))^{\ast}\frac{w^{(k)}_l}{\|w^{(k)}_l\|^2_2}
\frac{w^{(k)}_m}{\|w^{(k)}_m\|^2_2}w^{(k)}_n)=\]

\[=\sum_{r=0}^l\sum_{s=0}^m\mu(r,s,n;x^{-1},y^{-1})\tau^{\otimes_k}(\frac{w^{(k)}_{m+l+n-2(r+s)}}{\|w^{(k)}_{m+l+n-2(r+s)}\|^2_2}
\frac{w^{(k)}_l}{\|w^{(k)}_l\|^2_2} \frac{w^{(k)}_m}{\|w^{(k)}_m\|^2_2} w^{(k)}_n).\]
The coefficients depend only on $x,y$. The norms in the denumerator also coincide for all $k$. It only remains  to compute the trace.
As before, this reduces to linear combinations of powers of the initial trace (the canonical trace of the free group factor).
Since on group elements the trace is 0 or 1, we see that (\ref{eq23}) follows. Hence we have proved (\ref{eq1}). \rule{5pt}{5pt}

\section{\label{fn}The case of tensors of the Free Group Factors}

In what follows we will identify elements $g\in \mathbb{F}_N$ with their left regular representations $\lambda(g)\in \mathcal{L}(\mathbb{F}_N)$ and
we will use the notation from section \ref{Intro}.

The next technical lemma contains the calculations needed to prove the main result.
\begin{lemma} Let $k\ge 1$ and $w^{(k)}_n=\sum_{|v|=n}v^{\otimes_k}$. Let $x_1,...,x_{k},y_1,...,y_{k}$ be words in $\mathbb{F}_N$. Then:
\begin{equation}\label{ineq}
\sum_{n=0}^{\infty}\frac{1}{\|w^{(k)}_n\|^2_2}\| \mathbb{E}_k(x_1\otimes...\otimes x_{k}w^{(k)}_ny_1\otimes...\otimes y_{k})-
\end{equation}
\[ -\mathbb{E}_k(x_1\otimes...\otimes x_{k})\mathbb{E}_k(y_1\otimes...\otimes y_{k})w^{(k)}_n\|_2^2<\infty.
\]
\end{lemma}

{\it Proof: } We look first at the case $\mathbb{E}_k(x_1\otimes...\otimes x_k)\mathbb{E}_k(y_1\otimes...\otimes y_k)\ne 0$. By lemma \ref{nonzero},
it follows that $x_1=...=x_k=:x$ and $y_1=...=y_k=:y$. Let $l=|x|,m=|y|$. If $l=1$ or $m=1$ the inequality (\ref{ineq}) is trivial. Indeed,
suppose $l=1$ (the other case is even easier). Then we actually have
\[ \|\mathbb{E}_k(x^{\otimes_k}w^{(k)}_ny^{\otimes_k})-\mathbb{E}_k(x^{\otimes_k})\mathbb{E}_k(y^{\otimes_k})w^{(k)}_n\|^2_2=\]
\[\| \mathbb{E}_k(w^{(k)}_ny^{\otimes_k})-\mathbb{E}_k(y^{\otimes_k})w^{(k)}_n\|^2_2=\]
\[ \|w^{(k)}_n\mathbb{E}_k(y^{\otimes_k})-\mathbb{E}_k(y^{\otimes_k})w^{(k)}_n\|^2_2=0,\]
because $\mathcal{B}^{(k)}$ is abelian. We are left with the case $l,m\ge 2$. Using lemma \ref{k0}, the finiteness of the sum in (\ref{ineq})
follows in this case from the corresponding result in \cite{ss}.

For remainder of the proof assume  \[\mathbb{E}_k(x_1\otimes...\otimes x_k)\mathbb{E}_k(y_1\otimes...\otimes y_k)= 0.\] By the lemma \ref{nonzero},
there must be at least two distinct $x_{i_1},x_{i_2}$ and at least two distinct $y_{j_1},y_{j_2}$. The nonzero terms in the series (\ref{ineq})
\[\sum_{n=0}^{\infty}\frac{1}{\|w^{(k)}_n\|^2_2}\| \mathbb{E}_k(x_1\otimes...\otimes x_kw^{(k)}_ny_1\otimes...\otimes y_k)\|_2^2=\]

 \[\sum_{n=0}^{\infty}\frac{1}{\|w^{(k)}_n\|^2_2}\| \sum_{|v|=n}\mathbb{E}_k(x_1vy_1\otimes...\otimes x_kvy_k)\|_2^2\]
are given by the solutions $v$ of length $n\ge 0$  of the equations:
\[ x_1vy_1=x_2vy_2=...=x_kvy_k.\]
Observe that a solution exists only if $x_i=x_j\Leftrightarrow y_i=y_j$ for $1\le i,j\le k+1$. Without loss of generality,
assume $x_1\ne x_2(\Leftrightarrow y_1\ne y_2)$.
Any solution of length $n$ for the above system of equations is a solution of $v(y_1y^{-1}_2)=(x^{-1}_1x_2)v$.
Clearly $x^{-1}_1x_2\ne e\ne  y_1y^{-1}_2$. By lemma \ref{mylemma}, for each fixed
$n> {\rm max}\{|x_i| : 1\le i\le k+1\}+{\rm max}\{|y_j|:1\le j\le k+1\}=:n_0$ there is at most one solution to this equation. So
\[\sum_{n=n_0}^{\infty}\frac{1}{\|w^{(k)}_n\|^2_2}\| \sum_{|v|=n}\mathbb{E}_k(x_1vy_1\otimes...\otimes x_kvy_k)\|_2^2\le\]
\[\sum_{n=n_0}^{\infty}\frac{1}{\|w^{(k)}_n\|^2_2}\| \frac{w^{(k)}_n}{\|w^{(k)}_n\|^2_2}\|^2_2=\]
\[\sum_{n=n_0}^{\infty}\frac{1}{\|w^{(k)}_n\|^4_2} \stackrel{(\ref{norma wn})}{<}\infty.\qquad\rule{5pt}{5pt}\]

\begin{theorem}
For each $k\ge 1$,the conditional expectation $\mathbb{E}_k$ onto the Laplacian subalgebra $\mathcal{B}^{(k)}$ of
$\mathcal{L}(\mathbb{F}_N)^{\overline{\otimes}_k}$, generated by
\[ w^{(k)}_1=\sum_{v\in \mathbb{F}_N,|v|=1}v^{\otimes_k}, \]
is an asymptotic homomorphism and $\mathcal{B}^{(k)}$ is a strongly singular masa.
\end{theorem}
{\it Proof: } For $k=1$, $\mathcal{L}(\mathbb{F}_N)^{\overline{\otimes}_1}:=\mathcal{L}(\mathbb{F}_N)$ and the result was proved in
\cite{ss}. For $k\ge 2$, apply  theorem \ref{criteria}: use remark \ref{r1} and the proof for $k=1$ to verify the statement about the
existence of the $\ast$-isometry; use $Y=\{x_1\otimes ...\otimes x_{k}| x_i\in \mathbb{F}_N,1\le i\le k\}$ and $v_n:=w^{(k)}_n/\|w^{(k)}_n\|_2,n\ge 0$.
The theorem is a direct consequence of the above technical lemma and theorem \ref{criteria}.
\rule{5pt}{5pt}
\begin{acknowledgements} I want to thank my advisor, Florin R\u adulescu,  for his constant support, patience and encouragements.
I also want to thank Ionu{\c t} Chifan and Dorin Dutkay for valuable discussions and help.
\end{acknowledgements}


\begin{thebibliography}{20}
\bibitem{BR} Boca, F.; R\u adulescu, F. {\it Singularity of radial subalgebras in ${\rm II}_1$ factors asociated with free products of groups}, J. Funct. Anal., {\bf 103}  (1992), no. 1, 138--159.
\bibitem{Dix}Dixmier, J.  {\it Sous-anneaux amxmaux dans les facteurs de type fini}, Ann. Math., {\bf 59} (1954), 279--286.
\bibitem{Popa2} Popa, S. {\it Notes on Cartan subalgebras in type ${\rm II}_1$ factors}, Math. Scand., {\bf 57} (1985), 171--188.
\bibitem{Popa} Popa, S. {\it Singular maximal abelian $\ast$-subalgebras in continuous von Neumann algebras}, J. Funct. Anal., {\bf 103} (1983),  151--166.
\bibitem{Pytlik}, {\it Radial functions on free groups and a decomposition of the regular representation into irreducible components}, J. Reine Angew. Math., {\bf 326} (1981), 124--135.
\bibitem{FR}R\u adulescu, F. {\it Singularity of the radial subalgebra of $\mathcal{L}(\mathbb{F}_N)$ and the Puk\'ansky invariant},  Pacific J. Math. , {\bf 151}, no.2 (1991), pp 297-306.
\bibitem{GSS}Robertson, G.; Sinclair, A. M.; Smith, R. R. {\it Strong singularity for subalgebras of finite factors}, Intern. J. Math., {\bf 14}, No. 3 (2003) 235--258.
\bibitem{ss} Sinclair, A. M.; Smith, R. R. {\it The Laplacian MASA in a free group factor} Trans. Amer. Math. Soc. {\bf 355} (2003), no. 2, 465--475.
\bibitem{ss2} Sinclair, A. M.; Smith, R. R. {\it Strongly singular masas in type ${\rm II}_1$ factors} Geom. Funct. Anal. {\bf 12} (2002), no. 1, 199--216.
\end{thebibliography}
 \end{document}